\documentclass[10pt]{article}
\usepackage{amssymb,amsmath,amsfonts,bbm,pifont,upgreek,bbold,accents}  %

\usepackage[colorlinks=true]{hyperref}
\hypersetup{urlcolor=blue, citecolor=red}

\font\teneufm=eufm10
\font\seveneufm=eufm7
\font\fiveeufm=eufm5
\newfam\eufmfam
\textfont\eufmfam=\teneufm
\scriptfont\eufmfam=\seveneufm
\scriptscriptfont\eufmfam=\fiveeufm

\newcommand\beq[1]{ \begin{equation}\label{#1} }
\newcommand{\eeq}{ \end{equation} }
\newcommand{\beqno}{ \[ }
\newcommand{\eeqno}{ \] }
\newcommand\beqa[1]{ \begin{eqnarray} \label{#1}}
\newcommand{\eeqa}{ \end{eqnarray} }
\newcommand{\beqano}{ \begin{eqnarray*} }
\newcommand{\eeqano}{ \end{eqnarray*} }

\newtheorem{theorem}{Theorem}
\newtheorem{definition}{Definition}
\newtheorem{proposition}{Proposition}
\newtheorem{lemma}{Lemma}
\newtheorem{sublemma}{Sublemma}
\newtheorem{remark}{Remark}
\newtheorem{notationalremark}{Notational Remark}
\newtheorem{corollary}{Corollary}
\newtheorem{assumption}{Assumption}
\newtheorem{claim}{Claim}

\newtheorem{tools}{$\negsp\negsp$}[subsection]

\newcommand\thm[1]{ \begin{theorem}\label{#1}}
\newcommand\thmtwo[2]{ \begin{theorem}[#1]\label{#2}}
\newcommand\ethm{ \end{theorem} }
\newcommand\dfn[1]{ \begin{definition}\label{#1} \rm}
\newcommand\dfntwo[2]{ \begin{definition}[#1]\label{#2} \rm}
\newcommand\edfn{ \end{definition} }
\newcommand\pro[1]{ \begin{proposition}\label{#1}}
\newcommand\protwo[2]{ \begin{proposition}[#1]\label{#2}}
\newcommand\epro{ \end{proposition} }
\newcommand\lem[1]{ \begin{lemma}\label{#1}}
\newcommand\lemtwo[2]{ \begin{lemma}[#1]\label{#2}}
\newcommand\elem{ \end{lemma} }
\newcommand\sublem[1]{ \begin{sublemma}\label{#1}}
\newcommand\sublemtwo[2]{ \begin{sublemma}[#1]\label{#2}}
\newcommand\esublem{ \end{sublemma} }
\newcommand\rem[1]{ \begin{remark}\label{#1} \rm}
\newcommand\erem{ \end{remark} }
\newcommand\notrem[1]{ \begin{notationalremark}\label{#1} \rm}
\newcommand\enotrem{ \end{notationalremark} }
\newcommand\cor[1]{ \begin{corollary}\label{#1}}
\newcommand\cortwo[2]{ \begin{corollary}[#1]\label{#2}}
\newcommand\ecor{ \end{corollary} }
\newcommand\asmp[1]{ \begin{assumption}\label{#1}}
\newcommand\asmptwo[2]{ \begin{assumption}[#1]\label{#2}}
\newcommand\easmp{ \end{assumption} }
\newcommand\clm[1]{ \begin{claim}\label{#1}}
\newcommand\eclm{ \end{claim} }
\newcommand{\proof}{\par\medskip\noindent{\bf Proof\ }}
\newcommand\equ[1]{{\rm (\ref{#1})}}
\expandafter\chardef\csname pre amssym.def
at\endcsname=\the\catcode`\@
\catcode`\@=11
\def\undefine#1{\let#1\undefined}
\def\newsymbol#1#2#3#4#5{\let\next@\relax
 \ifnum#2=\@ne\let\next@\msafam@\else
 \ifnum#2=\tw@\let\next@\msbfam@\fi\fi
 \mathchardef#1="#3\next@#4#5}
\def\mathhexbox@#1#2#3{\relax
 \ifmmode\mathpalette{}{\m@th\mathchar"#1#2#3}%
 \else\leavevmode\hbox{$\m@th\mathchar"#1#2#3$}\fi}
\def\hexnumber@#1{\ifcase#1 0\or 1\or 2\or 3\or 4\or 5\or 6\or 7\or
8\or
 9\or A\or B\or C\or D\or E\or F\fi}
\ifcase\@ptsize
 \font\tenmsb=msbm10
 \font\sevenmsb=msbm7
 \font\fivemsb=msbm5
\or
 \font\tenmsb=msbm10 scaled \magstephalf
 \font\sevenmsb=msbm7 scaled \magstephalf
 \font\fivemsb=msbm5  scaled \magstephalf
\or
 \font\tenmsb=msbm10 scaled \magstep1
 \font\sevenmsb=msbm7 scaled \magstep1
 \font\fivemsb=msbm5 scaled \magstep1
\fi
\newfam\msbfam
\textfont\msbfam=\tenmsb
\scriptfont\msbfam=\sevenmsb
\scriptscriptfont\msbfam=\fivemsb
\edef\msbfam@{\hexnumber@\msbfam}
\def\Bbb#1{\fam\msbfam\relax#1}
\def\widehat#1{\setboxz@h{$\m@th#1$}%
 \ifdim\wdz@>\tw@ em\mathaccent"0\msbfam@5B{#1}%
 \else\mathaccent"0362{#1}\fi}
\def\widetilde#1{\setboxz@h{$\m@th#1$}%
 \ifdim\wdz@>\tw@ em\mathaccent"0\msbfam@5D{#1}%
 \else\mathaccent"0365{#1}\fi}

\def\RIfM@{\relax\ifmmode}
\def\nonmatherr@#1{\errmessage{\string#1\space allowed only in math mode}}
\def\Bbb{\RIfM@\expandafter\Bbb@\else
 \expandafter\nonmatherr@\expandafter\Bbb\fi}
\def\Bbb@#1{{\Bbb@@{#1}}}
\def\Bbb@@#1{\fam\msbfam\relax#1}
\def\setboxz@h{\setbox\z@\hbox}
\def\wdz@{\wd\z@}
\catcode`\@=\csname pre amssym.def at\endcsname
\newcommand{\giu}{{\medskip\noindent}}
\newcommand{\Giu}{{\bigskip\noindent}}
\newcommand{\nl}{{\smallskip\noindent}}
\newcommand{\noi}{{\noindent}}

\newcommand{\qed}{\hskip.5truecm
\vrule width 1.7truemm height 3.5truemm depth 0.truemm
\par\Giu}

\newcommand{\negsp}{\hspace{-.09truecm}}  

\newcommand{\dst}{\displaystyle}

\newcommand{\io}{{\infty }}

\newcommand\igl[2]{{ \int_{ {#1}}^{#2}  }}

\newcommand{\torus}{ {\Bbb T}   }
\renewcommand{\natural}{ {\Bbb N}   }
\newcommand{\real}{ {\Bbb R}   }
\newcommand{\integer}{ {\Bbb Z}   }
\newcommand{\complex}{ {\Bbb C}   }

\renewcommand{\a }{ {\alpha}   }
\renewcommand{\b}{ {\beta}   }
\newcommand{\g}{ {\gamma}   }

\renewcommand{\d}{ {\delta}   }

\newcommand{\e }{ {\varepsilon}   }

\renewcommand{\k}{ {\kappa}   }
\renewcommand{\l}{ {\lambda}   }
\renewcommand{\L}{ {\Lambda}   }
\newcommand{\m}{ {\mu}   }
\newcommand{\n}{ {\nu}   }
\newcommand{\x }{ {\xi}   }

\newcommand{\p}{ {\pi}   }

\renewcommand{\t}{ {\tau}   }
\newcommand{\f}{ {\varphi}   }

\renewcommand{\o}{ {\omega}   }
\renewcommand{\O}{ {\Omega}   }

\renewcommand{\Im}{{\, \rm Im\, }}

\newcommand{\cT}{ {\cal T} }

\newcommand{\cK}{ {\cal K} }

\newcommand{\cG}{ {\cal G} }
\newcommand{\cM}{ {\cal M} }

\newcommand{\cQ}{ {\cal Q} }

\newcommand\bks{\backslash}

\renewcommand\subset{\subseteq}

\newcommand\meas{{\,\rm meas\,}}
\newcommand\Lip{{\,\rm Lip}}

\newcommand\ttc{{\mathtt c}}
\newcommand\tty{{\mathtt y}}
\newcommand\ttx{{\mathtt x}}
\newcommand\ttH{{\mathtt H}}

\newcommand\ttP{{\mathtt P}}
\newcommand\ttB{{\mathtt B}}

\title{\bf
Nineteen Fifty-four: Kolmogorov's new `metrical approach' to Hamiltonian Dynamics
}

\begin{document}

\author{ 
L. Chierchia\footnote{Partially supporeted by the grant {\sl NRR-M4C2-I1.1-PRIN 2022-PE1-Stability in Hamiltonian dynamics and beyond-F53D23002730006-Finanziato dall'U.E.-NextGenerationEU}.} , I. Fascitiello
}

\maketitle

%

\begin{abstract}

\noindent
We review Kolmogorov's 1954 fundamental paper {\sl On the Conservation of Conditionally Periodic Motions under Small Perturbation of the Hamiltonian} (Dokl. akad. nauk SSSR,1954, vol. {\bf 98}, pp.527--530), both from the historical  and the mathematical point of view.
In particular, we discuss Theorem~2 (which deals with the  measure in phase space of persistent tori), the proof of which is not discussed at all by the author, notwithstanding its centrality in Kolmogorov's program in classical mechanics. 
\\
In Appendix, a recent interview to Ya. Sinai on KAM Theory is reported.
\end{abstract}

\nl

\nl

\vglue0.5truecm

\noi
Kolmogorov's 1954  paper {\sl On the Conservation of Conditionally Periodic Motions under Small Perturbation of the Hamiltonian} \cite{K} is probably one of -- if not, `the' -- most influential contribution to the modern development of classical mechanics and dynamical systems: In four pages, it started the celebrated KAM Theory, with precise statements and a clear outline of the main result. However, a complete discussion of this brief paper (four pages with a bibliography containing three items),  both from a historical and a mathematical point of view, 
is still missing\footnote{For a `friendly introduction' to the history and mathematics of KAM Theory, see \cite{dumas}.}. 

\nl
The plan of the present paper is the following.

\nl
In \S~1 ({\sl Historical remarks on Kolmogorov's influence on Classical Mechanics  in the 20th Century}), 
following  \cite{f}, and \cite{f22}, we try to put in a historical perspective Kolmogorov's revolutionary program in classical mechanics, as it emerges from his 1950's papers   and his lecture at the 1954 International Congress of Mathematicians in Amsterdam.

\nl
In \S~2 ({\sl The theorems in Kolmogorov's 1954 paper}),  we continue the mathematical discussion  started in \cite{C08}, considering, in particular,  Theorem~2, where Kolmogorov states that the Lebesgue measure of persistent invariant tori of analytic nearly--integrable Hamiltonian systems (with a non--degenerate integrable Hamiltonian and bounded phase space) tends to full measure, as the size of the perturbation goes to zero.
Theorem~2, unlike\footnote{Theorem~1  is followed by a precise outline of its proof, without, however,  including estimates and, in particular, without discussing the convergence of the Newton scheme; the missing analytical details have been discussed in 
 \cite{C08}, following closely Kolmogorov's outline.} Theorem~1, is not followed  by any mathematical details  or technical comments\footnote{In fact, after the statement of Theorem 2, Kolmogorov concludes his short paper simply with the following remark: 
 ``It seems that, in a sense, the `general case' is the case when the set $\cQ_\e$ [i.e., the set of all quasi--periodic trajectories, $\e$ being the size of the perturbation]  has an everywhere dense complement for all positive $\e$. 
Complications of this kind appearing in the theory of analytic dynamical systems were indicated in my paper \cite{K53} in connection with a more specific situation"}; we propose a proof of Theorem~2 based on the scheme of proof of Theorem~1 given in \cite{C08}.

\nl
Regarding Kolmogorov's second theorem (which, in some sense, was the real center of Kolmogorov's program),  the points of view of the two other founders of KAM theory appear to be rather different:

\nl
J.K.~Moser in his  1957 Mathematical Reviews report on Kolmogorov's paper \cite{K57},
says: 

\nl
\leftskip=1cm
\noindent
This very interesting theorem [i.e., Theorem~1 in \cite{K}] would imply that for an analytic canonical system which is close to an integrable one, all solutions but a set of small measure lie on invariant tori" (\underline{\href{https://mathscinet.ams.org/mathscinet/2006/mathscinet/search/publdoc.html?arg3=\&co4=AND\&co5=AND\&co6=AND\&co7=AND\&dr=all\&pg4=AUCN\&pg5=TI\&pg6=MR\&pg7=ALLF\&pg8=ET\&review\_format=html\&s4=\&s5=\&s6=MR0097598\&s7=\&s8=All\&sort=Newest\&vfpref=html\&yearRangeFirst=\&yearRangeSecond=\&yrop=eq\&r=1\&prev=t}{MR0097598}}). 

\leftskip=0cm

\nl
Moser does not make any reference here to Theorem~2, so it seems that he believed (in\footnote{Notice that the first contribution on small divisors by Moser is the famous 1962 paper \cite{Mo62} on area--preserving maps.} 1957) that Theorem 2 was a straightforward consequence of Theorem 1.

\nl
V.I. Arnol'd, on the other side, might have considered  the  lack of discussion of the proof of Theorem~2  the main motivation for  his celebrated paper entitled
{\sl `Proof of a theorem of A. N. Kolmogorov on the preservation of conditionally periodic motions under a small perturbation of the Hamiltonian'} \cite{A}, which included  a detailed discussion on the measure of the persistent tori, based, however,  on a  rather different scheme from that of Kolmogorov; in support of this comment, Sinai says\footnote{See the interview in Appendix.}:

\nl
\leftskip=1cm
\noindent
There were some gaps in the estimates of the measure of invariant sets [in \cite{K}]. That was the main point where Arnol'd complained about the proof by Kolmogorov. In Kolmogorov's paper, complete estimates of such as measure were not given.

\leftskip=0cm

\section{Historical remarks on Kolmogorov's influence on Classical Mechanics  in the 20th Century}

In the development of Mathematics, it is crucial to construct a narrative, as research progress is often articulated through reformulation of pre``existing ideas. Novel formulations frequently make the original concepts of past mathematicians unrecognizable, preventing access to meaning and intellectual engagement, which, in turn, provides cultural reference points.

\nl
It is with this objective that we historically analyze Kolmogorov's 1954 article \cite{K}, along with the text of the conference held by Kolmogorov in Amsterdam on September 9, 1954 \cite{K57}, within the broader historical context of the Soviet Union in the 1950s. We try to reconstruct the dissemination of the main ideas contained in both works and scrutinize their deeper significance through the testimonies of mathematicians close to Kolmogorov, such as Yakov G. Sinai and Vladimir I. Arnol'd, as well as mathematicians actively involved in the theory of dynamical systems, like Stephen Smale.

\nl
The International Congress of Mathematicians in Amsterdam, Netherlands, in 1954, ended with a lecture by Andrej Nikolaevich Kolmogorov (1903--1987). It was the second International Congress following the hiatus caused by political tensions and the Second World War, and the first to feature a Soviet delegation\footnote{Already suspended in 1936, the congress was reinstated only in 1950. However, on that occasion, the entire Russian academic community did not participate. In the proceedings of the ICM held in Cambridge, Massachusetts, within the Secretary's report section, the following is documented: 
\\
\noindent
Shortly before the opening of the Congress, the following cable was received from the President of the Soviet Academy of Sciences: ``The USSR Academy of Sciences appreciates having received a kind invitation for Soviet scientists to participate in the International Congress of Mathematicians to be held in Cambridge. Soviet mathematicians are very busy with their regular work, unable to attend the congress. I hope that the upcoming congress will be a significant event in mathematical science. Desire for success in congress activities. S. Vavilov, President, USSR Academy of Sciences." In \cite{proceedings50} p. 122.}. Notably, Kolmogorov's last overseas travels, dating back to the early 1930s\footnote{``en 1934 [...] quoique la fondation Rockefeller lui e\^ut accord\'e une bourse, Kolmogorov ne fut pas autoris\'e \`a se rendre \`a Paris pour travailler pr\'es d'Hadamard." In \cite{dem} p. 133.}, made this event historically significant. In Amsterdam, Kolmogorov presented a research program addressing open issues in classical mechanics within the broader framework of twentieth``century general theory of dynamical systems. The printed text \cite{K57} reads:

\nl
\leftskip=1cm
\noindent
My aim is to elucidate ways of applying basic concepts and results in the modern general metrical and spectral theory of dynamical systems to the study of conservative dynamical systems in classical mechanics. (p. 354).

\leftskip=0cm

\nl
This text comprises an extensive essay (approximately 24 pages), highly structured in argumentation, aimed at maintaining the reader's engagement. It is rich in bibliographic references to works published between 1917 and 1954, involving around 20 authors, primarily Soviet but also French and American.

\nl
The research program presented by Kolmogorov at the conference proved exceptionally fruitful in the latter half of the century. In the subsequent years in Moscow, Kolmogorov directed some of his students toward the field of dynamical systems, with also applications to celestial mechanics. With the fundamental contribution not only of Soviet mathematicians but also of others, this effort laid the foundation for what is now known as KAM theory.
The theory's name, formed from the acronyms of the names Kolmogorov, Arnol'd, and Moser, retrospectively reflects Kolmogorov's intent to establish a collective enterprise aimed at `young students of Moscow':

\nl
\leftskip=1cm
\noindent
My papers on classical mechanics appeared under the influence of von Neumann's papers on the spectral theory of dynamical systems and, particularly under the influence of the Bogolyubov--Krylov paper of 1937. 
\\
\noindent
[...] To accumulate specific information we organized a seminar on the study of individual examples. My ideas concerning this topic and closely related problems aroused wide response among young mathematicians in Moscow.  (In \cite{selworks} p. 521).

\leftskip=0cm


\nl
Lat us now focus on the historical reconstruction of the initial stage of the program, which includes Kolmogorov's papers between the spring of 1953 and the summer of 1954, months marked by upheavals in the Soviet Union following Stalin's death\footnote{March 5, 1953.}. The initial building blocks of the new mathematical landscape foreseen in the Amsterdam conference were published by Kolmogorov in two brief articles (4 pages each) \cite{K53} and \cite{K}, in the proceedings of the Soviet Academy of Sciences (Doklady Akademii Nauk SSSR), the first dated November 13, 1953, and the second -- arguably the most famous one -- on August 31, 1954, nine days before the Amsterdam conference. 
Both articles are referenced in the conference text, attesting to their significance in the research program. 
These short papers, almost devoid of bibliographical references, consist of the statements and hints at the proofs -- only for the first two (the third theorem is without a demonstration) -- of three theorems, without conceptual framing. Thus, we are faced with two distinct literary genres: the conference text, more literary and refined in its expressiveness, and concise twentieth--century research papers.


\nl
The second of Kolmogorov's two articles \cite{K} in the Proceedings of the Academy of Sciences, comprises two theorems. The first, which occupies the majority of the pages, is referred to as Kolmogorov's theorem on the persistence of invariant tori in Hamiltonian systems, following Arnol'd's nomenclature\footnote{{\sl Kolmogorov's 1954 theorem on the persistence of invariant tori under a small analytical perturbation of a fully integrable Hamiltonian system}. In \cite{A97}, p. 742.}.

\nl
Regarding the relevance in a  
longue dur\'ee
context of this theorem, Stephen Smale writes in \cite{smale}:

\nl
\leftskip=1cm
\noindent
It may be stated in conclusion that the outstanding unsolved problem in the ergodic theory is the question of the truth or falsity of metrical transitivity for general Hamiltonian systems. In other words the Quasi--Ergodic Hypothesis has been replaced by its modern version: the Hypothesis of Metrical Transitivity.
\\
This hypothesis played an important role in Birkhoff's later work. He not only believed it but part of his work is written assuming that it is true.
\\
$[...]$ These beliefs held sway in mathematical physics until Kolmogoroff's famous Amsterdam Congress paper in 1954 and subsequent work of Arnol'd and Moser in 1961--1962. The work of Kolmogoroff, Arnol'd, and Moser, KAM, showed that near `elliptic' closed orbits of a general Hamiltonian system on an energy surface, ergodicity failed. In that case there exist families of invariant tori of positive measure.
\\
Furthermore these elliptic orbits occur frequently in Hamiltonian systems. Thus the hypothesis of metrical transitivity is false in a definite way. (p.138--139).

\leftskip=0cm

\nl
The proof presented by Kolmogorov of his theorem on the persistence of invariant tori in this 1954 article \cite{K} has been examined \cite{C08}.\\
 However, within the same article \cite{K}, in its closing lines, there is a second theorem, neither elaborated nor demonstrated by Kolmogorov, which we will address in Section 2.

\nl
Both theorems, from a historiographical perspective, can be better understood in their cultural significance within the framework established by the research program presented at the Amsterdam conference. Conversely, this program gains greater concreteness when considering these initial steps taken independently by Kolmogorov.

\nl
It is in the conference text that we find the first explicit reference to the theorems, through a sentence attempting to condense their profound meaning within Kolmogorov's work:

\nl
\leftskip=1cm
\noindent
Theorems 1 and 2 in my paper [22]\footnote{He refers to \cite{K}.} assert that in the above--described situation the only change in the entire pattern for small $\theta$\footnote{For Kolmogorov's nomenclature, $\theta$ is the perturbation.} is that some of the tori corresponding to systems of frequencies for which the expression $(n, \lambda)$\footnote{$\lambda$ is the diophantine frequencies, $(n,\l)$ denotes inner product.} decreases too rapidly with increasing

$$|n|=\sqrt{\sum n^2_{\alpha}}$$

\noindent
may disappear while the majority of the tori $T^s_p$\footnote{$T^s_p$ are invariant unperturbed $s$--dimensional tori.}, retaining the character of motions occurring on them, are somewhat displaced in $\Omega^{2s}$\footnote{$\Omega^{2s}$ is the space phase.}, and still fill for small $\theta$ the region $G$\footnote{$G$ is a bounded region in $\Omega^{2s}$, where $H$ is defined.} to within a set of small measure. Thus, under small variations of $H$ the dynamical system remains non--transitive and the region $G$ continues to be decomposable, to within a residual set of small measure, into ergodic sets with discrete spectra (of the indicated specific nature). (\cite{K57}, p. 366).

\leftskip=0cm

\nl
Moreover, in the introduction of the text, the centrality of the Theorem 2 becomes evident:

\nl
\leftskip=1cm
\noindent
In conservative systems, asymptotically stable motions are impossible.
\\
\noindent
Therefore, for instance, the determination of individual periodic motions, however interesting it may be from the viewpoint of mathematics, has only a rather restricted real physical significance in the case of conservative systems. For conservative systems, the metrical approach\footnote{The title of our present paper derives from these words.} is of basic importance making it possible to study properties of a major part of motions. (\cite{K57}, p. 356).

\leftskip=0cm

\nl
The two works \cite{K} and \cite{K57} are therefore strongly interconnected: in the first (\cite{K}), we find the essential mathematics, distilled into two theorems, forming the foundation of Kolmogorov's research program; in the second (\cite{K57}), we delve into the detailed relevance of both theorems -- particularly, with more emphasis on the second theorem, which was not thoroughly analyzed in the paper \cite{K}-- within his research program. 

\nl
So, from a historiographical standpoint, it is crucial to analyze the dissemination of these contributions by Kolmogorov, which, in the initial years following their publication, primarily occurred through the distribution and translation of the Amsterdam conference text. 

\nl
The Proceedings of the 1954 ICM Congress in Amsterdam were published only in 1957 \cite{proceedings}. There, the original Russian text of Kolmogorov's conference is found, titled in Russian and French. In March 1958, Kolmogorov delivered a presentation at the Seminar on Analytical Mechanics and Celestial Mechanics hosted by Maurice Janet (1888--1983) \cite{K57-58} at the Faculty of Sciences of the Sorbonne, in Paris, on the same topic. The French translation of the Russian Kolmogorov's text \cite{K57} appears in the proceedings of the seminar of the young mathematician Jean--Paul Benz\'ecri (1932--2019). A note stated: \textit{L'auteur pr\'evoit la publication prochaine de d\`etails compl\`ementaires, dans un autre recueil}\footnote{p. 1, in a footnote. For Kolmogorov's trip to Paris, see to \cite{MAZ}.}, but this publication never materialized. In France, there was a certain interest in classical mechanics, yet it is noteworthy that all other seminars published in that issue of the Janet seminar, except for Kolmogorov's, dealt with general relativity\footnote{At his regard S. Dumas says: ``But KAM Theory [$\cdots$] also had the misfortune of playing out over roughly  the same interval during which the revolutions of modern physics took place.'' In \cite{dumas}, Preface, p. viii.}.

\nl
The original Russian version of Kolmogorov's conference was reprinted in the first volume of selected works edited by Vladimir Mikhailovich Tikhomirov (1934--) and published in 1985 by the publisher Nauka, two years before Kolmogorov's death.
Until 1972, there were only two circulating versions of the conference text: in Russian (1957) and in French (1958). In that year, the National Aeronautics and Space Administration (NASA) produced an English translation of the original Russian text. An English version, finally, was published within the English translation of the first volume of Kolmogorov's Selected Works, released in 1991 by Kluwer (the citation above is taken from this translation)\footnote{The author of this translation may be the Russian mathematical physicist Vladimir Markovich Volosov, who was in charge of the translation of the whole volume.}. The two articles in the Proceedings of the Academy of Sciences, also reprinted in the first volume of selected works in Russian, were translated for the English edition by Kluwer\footnote{The second one was published in English in 1977 in Stochastic Behavior in Classical and Quantum Hamiltonian Systems for the Volta Memorial conference, Como \cite{K.Como}.}.

\nl
The circulation of Kolmogorov's research program beyond the Iron Curtain was conditioned by the process of restoring contacts among mathematicians (where informal contacts resumed after a period of stagnation) and the resulting linguistic difficulties\footnote{The shift from French and German to Russian and English as the dominant languages within the international mathematical community after World War II created a new communication challenge. Starting in 1945, for instance, the British Mathematical Society had initiated the systematic English translation of the Russian journal ``Uspekhi Matematicheskikh Nauk", titled ``Russian Mathematical Surveys".}. Even more relevant, classical mechanics occupied a completely different r\^ole in many countries west of the Iron Curtain (except France) in the mid--century as compared to the prominent position it had held for a long time: equally neglected by theoretical physicists and mathematicians themselves, classical mechanics, as well as the study of dynamical systems, seemed confined to engineering issues\footnote{Clifford Truesdell in \cite{T} referred to this profoundly altered status of classical mechanics on multiple occasions. For instance: ``The word `classical' has two senses in scientific writing; (1) acknowledged as being of the first rank or authority, and (2) known, elementary, and exhausted (`trivial' in the root meaning of that word). In the twentieth century mechanics based upon the principles and concepts used up to 1900 acquired the adjective `classical' in its second and pejorative sense, largely because of the rise of quantum mechanics and relativity. [...]
Engineers still had to be taught classical mechanics, because in terms of it they could understand the machines with which they worked and could devise new machines for new purposes. Research in mechanics came to be slanted toward the needs of engineers and to be carried out largely by university teachers who regarded mathematics as a scullery--maid, not a goddess or even a mistress.'' (pp. 127--128).}. It was, in the end, a research program that implied a profound conceptual change intimately tied to the mathematical reformulation of classical mechanics -- it also had an impact on the three--body problem in celestial mechanics -- using new tools such as functional analysis and measure theory\footnote{In this regard, a historiographical issue arises concerning the origins of the conceptual shift or new paradigm proposed by Kolmogorov. For this, see \cite{f22}. Specifically, the ideas that Kolmogorov presented in 1954 are connected to research dating back to the 1930s. Precisely this span of approximately thirty years sheds light on the circumstances of mathematical research in the Soviet Union.}.

\nl
Let us return to the text of the paper \cite{K}.
The theorem concerning the persistence of invariant tori has given rise to contrasting interpretations, beginning with J\"urgen Moser's review published in ``Mathematical Reviews" in 1959\footnote{Moser reviews the conference text \cite{K57}, not the article \cite{K}. In reconstructing the paragraph structure and topics covered, he observes that ``at the heart of this discourse lies the author's novel assertion concerning the conservation of conditionally periodic solutions," and concludes the review with, ``The proof of this theorem was published in Dokl. Akad. Nauk SSSR 98 (1954), 527--530 [MR0068687], yet the discussion on convergence appears unconvincing to the reviewer. This very interesting theorem would imply that for an analytic canonical system which is close to an integrable one, all solutions but a set of small measure lie on invariant tori".}.
These doubts revolve around the validity of the proof of the theorem on the persistence of invariant tori -- the first of the theorems contained in the paper -- and the true original statement. Indeed, frequent reference is made to a `KAM theorem', both in relation to Moser and Arnol'd himself. Arnol'd, in the passage where he uses the name we adopt here, is unequivocal:

\nl
\leftskip=1cm
\noindent
This theory is referred to as KAM, or Kolmogorov--Arnol'd--Moser, and it is commonly asserted that there exists a KAM theorem. I have never managed to discern the specific theorem in question. In \cite{A04}, p. 622. 

\leftskip=0cm

\nl
Further insight regarding Kolmogorov's demonstration of this theorem can be found in \cite{sinai}, an article authored by Sinai within the volume \textit{Kolmogorov in Perspective}, a compilation of works penned by colleagues and students of Kolmogorov concerning his contributions and life:

\nl
\leftskip=1cm
\noindent
In the fall of 1957 I became a graduate student under Andrei Nikolaevich. At the same time he began a famous course of lectures on the theory of dynamical systems, which later was continued as a seminar. Much has already been written about this seminar. Among those present, besides us, were V. M. Alekseev, V. I. Arnol'd, L. D. Meshalkin, M. S. Pinsker, M. M. Postnikov, K. A. Sitnikov, and many others. The first part of the course definitely had a probabilistic bias, although in presenting the von Neumann theory of dynamical systems with pure point spectrum Kolmogorov made use of Pontryagin's theory of characters. For probabilists these were completely unfamiliar objects, of course, but he used them as freely as everything else. Before beginning the lecture he asked the listeners who was familiar with the theory of characters, and only Postnikov and Sitnikov raised their hands. Later in the course he presented the theorem that was to become the basis for the famous KAM theory, together with a complete proof. In early 1958 Andrei Nikolaevich departed to spend half a year in Prance and left Meshalkin and me a program for preparation for the examination in classical mechanics, which included this proof. (p. 117--118).

\leftskip=0cm

\nl
The origin of Moser's and others' `doubts' can be partly attributed to the fact that the original sources became available in English many years after the original publication. However, it seems that interpretative oscillations can be rather linked to two profound conceptual circumstances:

\begin{enumerate}
    \item The first, more common and interesting, is related to the level of details required for a proof to be truly convincing and the limits of the general conception of absolute deductive proof. In \cite{C08}, the proof of this theorem proposed in the original paper has been reconstructed by completing the missing details.
    \item The second circumstance refers to the understanding and interpretation of the theorem's statement itself, in connection with subsequent reformulations or theorems inspired by Kolmogorov's. Indeed, the statement of the theorem, along with the second theorem found in \cite{K}, implicitly contains the core of Kolmogorov's research program. However, it has often been interpreted without reference to the dense presentation of this program contained in the conference text (also because this was less accessible). Since article \cite{K} encompasses revolutionary features, a better understanding of them could be gained by considering Kolmogorov's detailed presentation accompanied by bibliographic references.
\end{enumerate}

\noindent
In reference to this second aspect, it appears useful to reconsider the second theorem contained in Kolmogorov's second article \cite{K}, and in particular its proof; see  \S~\ref{theorem2} below.

\section{The theorems in Kolmogorov's 1954 paper}

Here, we analyze in detail the two theorems appearing in Kolmogorov's paper \cite{K}: in \S~\ref{theorem1}
we recall the discussion made in \cite{C08} of Theorem~1, adding a few remarks; in \S~\ref{theorem2} we propose a proof of Theorem~2 based on the scheme of proof of Theorem~1 given by Kolmogorov and implemented in \cite{C08}.

\subsection{Theorem~1}\label{theorem1}

Theorem~1 in \cite{K}  is the celebrated Kolmogorov's  Theorem on the persistence of invariant tori for analytic Hamiltonian systems, from which stemmed 
KAM Theory. The following is an extended statement based on a proof (\cite{C08}), which completed the missing analytical details  along the original outline.
Such theorem deals with small analytic perturbations of a real analytic Hamiltonian in `Kolmogorov normal form', namely, a Hamiltonian $K$ of the form\footnote{$\o\cdot y=\sum_j \o_j y_j$ is the standard inner product and $Q=O(|y|^2)$ means that $Q$ vanishes together with its $y$--derivatives at $y=0$.}
\beq{K}
K=K(y,x)=E+\o\cdot y +Q(y,x)\,,\qquad Q=O(|y|^2)\,,
\eeq
where $E\in\real$,  
\beq{omega}
\o\in\real^d_{\g,\t}:= \big\{\o\in\real^d:\ |\o\cdot n|\ge \frac\g{|n|^\t}\,, \quad \forall\ n\in\integer^d\bks\{0\big\}\,,
\eeq is a  Diophantine  frequency, for some $\t\ge n-1$, $\g>0$, and  $Q$  is non--degenerate in the sense that 
\beq{Knd}
\det \langle \partial^2_y Q(0,\cdot)\rangle:=\igl{\torus^d}{} \partial^2_y Q(0,x) \frac{dx}{(2\p)^d}\neq 0\,.
\eeq
The phase space is\footnote{$B_\x(y)$ denotes the Euclidean $d$--ball with radius $\x$, centred at $y_0$ and $\torus^d:=\real^d/(2\pi \integer^d)$ is the standard flat $d$--dimensional torus.} $\cM=B_\x(0)\times \torus^d$, endowed with the standard symplectic form $dy\wedge dx =\sum_j dy_j\wedge dx_j$, which means that the Hamiltonian flow generated by the Hamiltonian $K$, $t\to \Phi_K^t(y,x)$, is the solution of the system
so that the solution of  the Cauchy problem 
\beqno
\left\{
\begin{array}l
\dot y(t) = - \partial_x K(y(t),x(t);\e)\,,\\
\dot x(t) =  \partial_y K(y(t),x(t);\e)\,,
\end{array}\right.
\ 
\left\{
\begin{array}l
y(0) =y\,,\\
x(0) =  x\,.
\end{array}\right.
\eeqno
The special feature of a Hamiltonian $K$ in Kolmogorov's normal form is that  {\sl the torus $\cT_0:=\{0\}\times \torus^d\subset \cM$ is a Lagrangian transitive invariant  torus for  $K$}, since $\Phi_K^t(0,x)= (0,x+\o t)$. The Diophantine vector $\o$ is called {\sl the frequency vector of the invariant torus $\cT_0$}.

\nl
Given $\x,\e_0>0$,   define the  complex domain
\beq{W}
W_{\x,\e_0}:= D^d_\x(0)\times \torus^d_\x\times D^1_{\e_0}(0)\subset \complex^{2d+1}\,,
\eeq
where $D^m_r(z)$ denotes the complex $d$--ball of radius $r$ centered at $z\in \complex^m$, and $\torus^d_\x$ is the complex neighbourhood of the torus $\torus^d$ given by $\{x\in \complex^d: |\Im x_j|<\x\,, \forall j\}/(2\p \real^d)$. 
For a real analytic function $f:W_{\x,\e_0}\to \complex$ we denote its sup norm on $W_{\x,\e_0}$ by $\|f\|_{\x,\e_0}$, and its sup--norm (at fixed $\e$) by $\|f\|_\x$. Then, Theorem~1 in \cite{K} can be formulated as follows\footnote{Part (i) is essentially Kolmogorov's original statement, part (ii) contains the associated estimates.}.

\giu
{\bf Theorem 1} (i) {\sl Let $\o\in\real^d_{\g,\t}$ and $K$ be a Hamiltonian in Kolmogorov's normal form as in \equ{K}--\equ{Knd} with $K$ real analytic and bounded on $W_{\x,\e_0}$ for some $\x,\e_0>0$; let $P=P(y,x;\e)$ be a real analytic function   on $W_{\x,\e_0}$. Then, for any $0<\x_*<\x$, there exists  $0<\e_*\le \e_0$ and, for any $0\le \e<\e_*$, a  near--to--identity symplectic transformation $\phi_*:D_{\x_*}^d(0)\times\torus^d_{\x_*}\to D^d_\x(0) \times \torus^d_\x$, real analytic on $W_{\x_*,\e_*}$, 
such that  the Hamiltonian $H\circ \phi_*$, where 
$H:=(K+\e P)$, is in Kolmogorov normal form:
\beq{Kstar}
H\circ \phi_*= K_*=E_*+\o\cdot y+Q_*\,,\qquad Q_*=O(|y|^2)\,.
\eeq
{\rm (ii)}  In the above statement one can take  $\e_*=\min\{\e_0,\ttc_*^{-1}\}$  where
\beq{estar}
\ttc_*= \ttc \g^{-4} (\x-\x_*)^{-\n} C^\n\, \|P\|_{\x,\e_0}\,,\quad {\rm and}\quad  
\left\{\begin{array}{l}
C:=\max\big\{|E|, |\o|,\|Q\|_{\x,\e_0}, \|T\|\,,1\big\}\,,\\
T:=\langle \partial^2_y Q(0,\cdot)\rangle^{-1}\|\,,
\end{array}\right.
\eeq
and $\ttc,\n>1$  are suitable constants depending only on $d$ and $\t$.
 Furthermore, for any complex $\e$ with $|\e|<\e_*$, one has 
\beqno
\|\phi_*-{\rm id}\|_{\x_*}\,, |E-E_*|\,,\|Q_*-Q\|_{\x_*}\,,  \|\langle \partial^2_y Q(0,\cdot)\rangle^{-1}- \langle\partial^2_y Q_*(0,\cdot)\rangle^{-1}\| \le   \ttc_*\, |\e|\,.
\eeqno
}

\nl
Let us make a few remarks.

\giu
{\bf (1.1)} ({\sl On the dependence of the smallness condition upon the Diophantine constant $\g$})\\
A detailed proof, {\sl apart} form the explicit dependence upon the Diophantine constant $\g$ (which plays an important r\^ole in 
the analysis of the measure of persistent tori), based on Kolmogorov's original outline, has been given in \cite{C08}: compare, in particular,   Lemma~5,  
Eq. (27) (the factor 2 in the definition of $C$ in Eq. (26) has, here, been absorbed in the constant $\ttc$), and Eq. (31). The way the constant $\ttc_*$ depends upon $\g$ needs a short discussion.

\nl
The dependence upon $\g$ comes in through the constant $\bar c$ in
Eq.~(18) of \cite{C08} (beware that the Diophantine constant $\g$ is denoted $\k$ in \cite{C08}). Now, in the first line of Eq.~(18) one can actually take $\bar c=\g^{-2} \bar c_0$ with $\bar c_0=\bar c_0(d,\t)$ depending only on $d$ and $\t$, since the factor 
$\g^{-1}$ appears every time the small--divisor operator  $D^{-1}_\o$ (i.e., the inverse of the directional derivative $D_\o=\sum_j \o_j \partial_{x_j}$ acting on zero--average, real analytic functions on $\torus^d$) is applied, and  the formulae defining the functions in the left hand side of (18) involve $D^{-1}_\o$ at most {\sl twice}; compare the formulae at the beginning of p. 135 of \cite{C08}. Then, it is easy to check that the constant $\bar c$ in the estimate on the norm of  the  `new' perturbing function $P'$ in the second line of Eq.~(18) can be taken to be\footnote{The factor $(\g^{-2})^2$ comes from the term $P^{(1)}$; compare Eq.~(11).} $\bar c= \g^{-4} \bar c_1$, with $\bar c_1=\bar c_1(d,\t)$. Therefore also, $c$ in Eq.~(22) and $c_*$ in Eq.~(27) in \cite{C08} are proportional to a constant $\bar c_*(d,\t) \g^{-4}$, which leads to \equ{estar} above.\\
Incidentally, we observe  that the argument sketched here  shows that
{\sl  the relation $c \e_* \g^{-4}\|P\|_\x<1$, with a constant $c$ independent of $\g$,  cannot be improved following Kolmogorov's scheme}: Indeed,   
the norm of $\|P'\|$ cannot be estimated better than by  $\g^{-4}\|P\|^2$ times a constant independent of $\g$, and  iterating this relations (i.e., replacing $\|P\|$ with $\|P_{j-1}\|$, $\|P'\|$ with $\|P_j\|$; $P_0:=P$), one finds that $|\e^{2^j}|\|P_j\|\sim \g^4 (|\e|\g^{-4}\|P\|)^{2^j}$, so that, in order for the Newton scheme to converge, it is {\sl necessary} that $c\,\e_* \g^{-4}\|M\|<1$. \\
On the other hand, following Arnol'd's approach \cite{A} -- which is a Newton scheme based on approximate solutions of Hamilton--Jacobi equations, where the new perturbing function is of order $\e^2 \g^{-2} \|P\|^2$ --  allows for a final  condition of the form $c\, \e_* \g^{-2}<1$, which turns put to be optimal (as far as {\sl primary tori} are concerned\footnote{Primary tori are invariant tori which are a deformation of integrable tori and which, in particular, are graphs over $\torus^d$; for a discussion of a KAM Theory for {\sl primary and secondary} tori, see \cite{BC23}.});  compare, e.g., \cite{CK1}.

\Giu
{\bf (1.2)}  ({\sl On the structure of Kolmogorov's transformation})

\nl
Kolmogorov's transformation $\phi_*$ has a particularly simple form. Indeed, Kolmogorov describes in detail the transformation $\phi_1$, which is the first transformation of the iteration, conjugating the starting Hamiltonian $H=K+\e P$ to a new Hamiltonian $H_1:=K_1+\e^2 P_1:=H\circ\phi_1$, with $K_1$ in Kolmogorov's normal form  with same\footnote{After the description of $\phi_1$ Kolmogorov adds (\cite[p. 55]{K}): ``The construction of further approximations is not associated with new difficulties. Only the use of condition \equ{omega} for proving the convergence of the recursions, $\phi_j$, to the analytical limit for the recursion $\phi_*$ is somewhat more subtle."
} $\o$.
\nl
Now, the transformation $\phi_1$ belongs to the (formal) group of near--to--identity symplectic transformations $\cG$ of the form
\beqno
\phi:(y',x')\mapsto\left\{
\begin{array}{l} y= y'+\e \big(u(x')+U(x')  y'\big)\\
x=x'+\e \a(x')\
\end{array}
\right.
\eeqno
with $U$  a $(d\times d)$  matrix (depending periodically on $x'$): such transformations are defined, 
for small $\e$, in a neighborhood of the origin   times $\torus^d$; compare Remark~2, and in particular, Eq.~(9),  in \cite{C08}. In the recursion,
$\phi_j$ will have the same form but with $\e$ replaced by $\e^{2^{j-1}}$, and\footnote{Of course, all the symbols indexed by $*$ depend on $\e$ (and on the fixed $\o$).} $\phi_*=\lim_j\phi_1\circ\cdots \phi_j$ will be given by 
$$(I,\theta)\mapsto \phi_*(I,\theta) =(I,\theta)+ \e\big(u_*(\theta)+U_*(\theta) I, \theta+\e \a_*(\theta) \big)\in \cG\,.$$ 
Thus, defining 
\beq{zeta}
\zeta_*(\theta):= \phi_*(0,\theta)=\big(\e u_*(\theta),\theta+\e \a_*(\theta) \big)\,,
\eeq 
the final invariant torus for the original Hamiltonian $H$
is given by 
\beqno
\cT_*:=\big\{(y,x)= \zeta_*(\theta):  \theta\in \torus^d \big\}\,,\qquad {\rm and}\qquad \Phi_H^t\big(\zeta_*(\theta)\big)= \zeta_*(\theta+\o t)\,.
\eeqno
Observe that, since the map $\theta\mapsto \theta+\e\a_*(\theta)$ is a diffeomorphism of $\torus^d$ with inverse of the form $x\mapsto x+\e a_*(x)$, the 
{\sl invariant torus $\cT_*$ is a graph over $\torus^d$ given by 
\beqno
\cT_*=\big\{(y,x)=\big(\e \bar y_*(x),x\big): x\in \torus^d\big\}\,,\qquad \bar y_*(x):=u_*(x+\e a_*(x)) \,.
\eeqno
}

\noi
{\bf (1.3)}  ({\sl On  $\e$--analyticity and the convergence of Lindstedt series})

\nl
Let us make the obvious remark  -- which, however, seems to have been completely overlooked! -- that from Theorem~1, it follows immediately that {\sl the invariant torus $\cT_*$ depends analytically on $\e$}, since $\phi_*$ is real analytic on $W_{\x_*,\e_*}$, as the above function $\zeta_*$ is  analytic in $\{\e\in \complex: |\e|<\e_*\}$.

\nl
This observation implies at once that {\sl the Lindstedt series proposed for the first time in} \cite{Lin} -- i.e., the formal $\e$--expansion of quasi--periodic trajectories for nearly--integrable Hamiltonian systems (which in the present setting is given by $\zeta_*$) -- {\sl are actually convergent $\e$--power series}, a fact that was formally settled, after eighty years from Lindstedt's memoirs and thirteen years after Kolmogorov's paper, by J. Moser in 1967 \cite{M67} using his version of KAM theory (which, again, is rather different from Kolmogorov's approach).
\\
Incidentally, it is worthwhile to mentioned that H.~Poincar\'e, apparently, thought that Lindstedt series were divergent, as it appears from his comments in \cite{Poi}:
``M. Lindstedt ne d\'emontrait pas la convergence des d\'eveloppements qu'il avait ainsi form\`es, et, en effet, ils sont divergents" (\cite{Poi}, vol. II, \S~IX, n. 123); and later: 
``Il semble donc permis de conclure que le s\'eries (2) ne convegent pas. Toutfois le raisonnement qui pr\'ec\`ede ne
suffit pas pour \'etablir ce point avec une rigueur compl\`ete[...] Tout ce qu'il m'est permis de dire, c'est qu'il est fort invraisemblable". (\cite{Poi}, vol. II, \S~XIII entitled `Divergence des series de M. Lindstedt', n. 149).

\subsection{Theorem 2}\label{theorem2}
In Theorem~2,     Kolmogorov considers real analytic nearly--integrable Hamiltonian systems, namely, one--parameter families of Hamiltonian systems governed
 by a real analytic Hamiltonian
\beq{H}
(\tty,\ttx,\e)\in W:=V\times\torus^d\times (-\e_0,\e_0)\mapsto 
\ttH(\tty,\ttx;\e):=\ttH_0(\tty)+\e \ttP(\tty,\ttx;\e)\,,
\eeq
where  $V\subset \real^d$ is a bounded regular open connected set, and  $\e_0>0$; `regular', here, means that\footnote{
$\meas$ denotes Lebesgue measure.
This regularity assumption on the (boundary of) the set $V$ is not present in Kolmogorov's paper: Kolmogorov speaks simply of  `bounded regions'.} 
\beqno
\lim_{\d\to 0}\meas(V\bks V^{(\d)})=0\,,\qquad {\rm where}\qquad V^{(\d)}:= \{\tty\in V: B_\d(\tty)\subset V\}\,.
\eeqno
The phase space is the set  $\cM:=V\times \torus^d$, endowed with the standard symplectic form $d\tty\wedge d\ttx=\sum_j d\tty_j\wedge d\ttx_j$, and $\e$ is a small parameter. Denote by
$\phi_\ttH^t(\tty,\ttx)$ the Hamiltonian flow starting at $(\tty,\ttx)\in \cM$.\\
In considering such systems, Kolmogorov says: 

\nl
``There arises the natural hypothesis that at small $\e$ the `perturbed tori' obtained by Theorem~1 fill the larger part of the region $\cM$. This is also confirmed by Theorem~2, pointed out later". 

\nl
Then, Kolmogorov defines 
the set $\cQ_\e$ of Hamiltonian trajectories  in $\cM$, which are quasi--periodic with  frequencies $\o\in\real^d$, 
i.e., trajectories of the form $\phi_\ttH^t(\tty,\ttx) =(Y(\o t), X(\o t))$ for suitable analytic functions $\theta \in \torus^d\mapsto (Y(\theta),X(\theta))\in \cM$, and,
at the end of \cite{K}, states the following

\Giu
{\bf Theorem 2}{\sl\  Let $\ttH$ be as in \equ{H} and assume  $\det \partial_\tty^2 \ttH_0\neq 0$ on $V$. Then, $\dst\lim_{\e\to0}\meas (\cM\bks \cQ_\e)=0$.
}

\giu
As already mentioned, this statement is not accompanied by any remark, nor references. 
In the rest of this section, we will show how one can deduce Theorem~2 from Theorem~1 and its proof. 

\Giu
\proof {\bf of Theorem 2}

\nl
{\bf (2.1)} {\sl Local reduction} 

\nl
The claim of Theorem~2 is actually of local nature. Indeed, 
since $V$ is a regular set, it is enough to show that, for each $\d>0$, $\lim_{\e\to 0} \meas((V^{(\d)}\times \torus^d)\bks \cQ_\e)=0$. 
Furthermore, since $\ttH$ is real analytic on $V\times \torus^d$ and $V^{(\d)}$  is compact, $\ttH$ is real analytic and bounded on $\cup_{\tty\in V^{(\d)}} D^d_{\x_0}(\tty)\times \torus^d_{\x_0}$ for a suitable $0<\x_0<\d$ (for all $|\e|<\e_0$). Also, since
$\det \ttH_0''\neq 0$ on $V$, by the Implicit Function Theorem, there exists $0<r<\x_0/2$ such that the unperturbed frequency map 
\beqno
\tty\in V\mapsto \o_0(\tty):= \partial_\tty \ttH_0(\tty)\,,
\eeqno
is an analytic diffeomorphism from $B$ onto $\O:=\o_0(B)$, for any closed ball $B=B_r(\tty)$ with $\tty\in V^{(\d)}$. Therefore,  since $V^{(\d)}$ can be covered by a finite number of such balls $B$, {\sl it is enough to prove that $\lim_{\e\to 0} \meas ((B\times \torus^d)\bks \cQ_\e)=0$, for  any set $B=B_r(\tty)$ with $\tty\in V^{(\d)}$}.

\giu
{\bf (2.2)} {\sl Application of Theorem~1}

\nl
At \cite[p. 55]{K}, Kolmogorov says: ``The condition of the absence of `small denominators (3) [i.e., the Diophantine inequalities in  \equ{omega}] should be considered, `generally speaking', as fulfilled since for any $\t>d-1$ for all points of a $d$--dimensional space $\o=(\o_1,...,\o_d)$ except the set of Lebesgue measure zero it is possible to find $\g=\g(\o)$ for which $|\o\cdot n|\ge \g/|n|^\t$ whatever the integers $n\neq 0$ are".
\\
Indeed,  it is an elementary observation that, if $\t>d-1$ and we define $\O_\g:=\{\o\in \O: \o\in \real^d_{\g,\t}\}$, then 
\beq{O}
\meas (\O\bks \O_\g)\le c \g\,,
\eeq
for a suitable constant $c$ depending on $d, \t$ and on the diameter of\footnote{If $\d_\O$ denotes the diameter of $\O$, one has $\O\bks\O_\g\subset \big\{\o\in \O: \exists n\neq 0\, {\rm s.t.}\, \big|\o\cdot \frac{n}{|n|}\big|< \frac\g{|n|^{\t+1}}\big\}$, which implies $\dst \meas (\O\bks \O_\g)\le \sum_{n\neq 0} \frac{\g}{|n|^{\t+1}} \d_\O^{d-1}=:c \g$.
} $\O$. 

\nl
To proceed in the discussion, fix $\d>0$ and pick a ball $B=B_r(\tty)$ as in {\bf (2.1}) above; fix
(once and for all) $\t>d-1$ and let $\g>0$ (eventually, $\g$ will be chosen as a suitable power of $\e$). Denote
\beqno
\ttB_\g:=\big\{\tty\in B: |\o_0(\tty)\cdot n|\ge  \frac\g{|n|^\t}\,, \ \forall\ n\in\integer^d\bks\{0\}\big\}\,, \ 
\O_\g:=\o_0(\ttB_\g)=\{\o\in \o_0(B): \o\in \real^d_{\g,\t}\}\,.
\eeqno
Observe that $B_\g$ and $\O_\g$ are nowhere dense  sets and that $\o_0$ is a Lipeomorphism (bi--Lipschitz homeomorphism) between them, being an analytic diffeomorphism of $B$ onto $\O=\o_0(B)$.
Fix $\tty\in \ttB_\g$ and consider the trivial symplectic map 
\beqno 
\phi_0:(y,x)\in B_\x(0)\times \torus^d \mapsto \phi_0(y,x;\tty):=(\tty+y,x)\,,
\eeqno
where $\x:=\x_0/2$.
Then,  define
\beqno
H(y,x):= \ttH_0\circ \phi_0+\e \ttP\circ \phi_0
=:K+\e P\,,\quad P=P(y,x;\e):=\ttP(\tty+y,x;\e)\,,
\eeqno
and observe that $H$ is real analytic and bounded on\footnote{Recall the definition in \equ{W}, and that, by {\bf (2.1)}, $\ttH$  is real analytic and bounded on $\cup_{\tty\in V^{(\d)}} D^d_{\x_0}(\tty)\times \torus^d_{\x_0}$ for all $|\e|<\e_0$.} $W_{\x,\e_0}$.
By Taylor's formula, one has
\beqno
K:= E+\o\cdot y+Q\,,\quad {\rm with}\quad \left\{
\begin{array}{l}
 E:=\ttH_0(\tty)\,,\quad\o:=\o_0(\tty) \\
\dst  Q=\Big(\igl{0}1 (1-t)\ttH''_0(\tty+t y)dt\Big)y\cdot y\,.
 \end{array}
 \right.
\eeqno
Since for any $\tty \in \ttB_\g$, $\o\in\O_\g$, {\sl we can apply Theorem~1} to $H$ and   get 
a near--to--identity symplectic transformation $\phi_*$   so that \equ{Kstar} holds for any  $\e<\e_*$. Notice that {\sl everything here} ($H$, $\phi_*$, etc.) 
{\sl is parameterized by} $\tty\in \ttB_\g$.
Thus, 
$$\cT_* =\cT_*(\tty):=\psi(\{0\}\times\torus^d)\,,\qquad{\rm where}\quad  \psi:=  \phi_0\circ \phi_*\,,$$ 
is a real analytic Lagrangian torus invariant for the flow of $\ttH$ and spanned by  Diophantine quasi--periodic trajectories. In fact, defining the `Kolmogorov's transformation'
\beq{Kt}
\psi_{\!{}_K}: (\tty,\theta)\in \ttB_\g\times\torus \mapsto \psi_{\!{}_K}(\tty,\theta):=  \psi(0,\theta;\tty)\stackrel{\equ{zeta}}{=} \big(\tty+\e u_*(\theta;\tty), \theta+\e\a_*(\theta;\tty)\big)\,,
\eeq
we find 
\beq{Kflow}
t\mapsto \Phi^t_\ttH\psi_{\!{}_K}(\tty,\theta) = \psi_{\!{}_K}(\tty,\theta+\o t)\,.
\eeq

\noi
{\bf (2.3)} {\sl The Kolmogorov's set}

\nl
In view of \equ{O}, in  order to get a full measure set as $\e$ goes to zero, it is natural to choose $\g$ as a suitable power of $\e$  so that the smallness condition of Theorem~1 holds {\sl uniformly} in phase space. For example,  
 if we take $\g=\e^{1/5}$, we see that $\e_*$ in point (ii) of Theorem~1,  for $\e$ small enough,  is given by 
$\e_*\sim \e^{4/5}$, so that the condition  $\e<\e_*$  is fulfilled for any $\e>0$ small enough and any $\tty\in \ttB_\g=\ttB_{\e^{1/5}}$. 
With these choices, the set
\beq{Kset}
\cK_\e:= \psi_{\!{}_K}(\ttB_\g\times\torus^d)\,,\qquad \g=\e^{1/5}\,,
\eeq
defines a set of invariant tori for $\ttH$, which, by \equ{Kflow}, is made up of quasi--periodic trajectories, so that $\cK_\e\subset \cQ_\e$.
Also, from \equ{O}, since $\ttB_\g=\o_0^{-1}(\O_\g)$ and $\o_0$ is a diffeomorphism, it follows that
\beq{measB}
\meas \big((\ttB\times\torus^d) \bks (\ttB_\g\times\torus^d)\big)\le c'\ \e^{1/5}\,.
\eeq
All this, in our opinion, must have been rather obvious to Kolmogorov. Furthermore, the Kolmogorov's map $\psi_{\!{}_K}$ in \equ{Kt} is a near--to--identity map and it is very tempting, at this point, to conclude that also $\meas ((\ttB\times\torus^d) \bks \cK_\e)\to 0$ as $\e\to 0$ concluding the proof of Theorem~2. 
Clearly, to complete the argument, one needs to have more information on the regularity of $\psi_{\!{}_K}$ in order to control how Lebesgue measure changes under its action. For example, one can  check that  $\psi_{\!{}_K}$ is a lipeomorphism with Lipschitz constant  arbitrarily close to 1 as $\e\to 0$, in which case, from \equ{Kset} and \equ{measB},  Theorem~2 follows immediately.

\giu
{\bf (2.4)} {\sl Lipschitz properties}

\nl
What one needs to show is that the functions $u_*$ and $\a_*$ are Lipschitz functions with uniformly bounded Lipschitz constants on $\ttB_\g$. Now, the way $\tty$ enters in the construction of $\phi_*$ is only through $\o =\o_0(\tty)$, and, since $\o_0$ is an analytic function, it is enough to check that  
{\sl $\phi_*$ (and hence $u_*$ and $\a_*$) is a Lipschitz functions of $\o$ with uniformly bounded Lipschitz constants on $\O_\g$.} 

\nl
The starting simple  observation is that if $u=\sum_{n\neq 0} u_n e^{i n\cdot x}$ is an analytic map on $\torus^d$ with zero average, then
\beq{Dom}
(D_\o^{-1} u) (x):= \sum_{n\neq 0} \frac{u_n}{i \o\cdot n} e^{i n\cdot x}\,,
\eeq 
depends in a Lipschitz way on $\o\in \O_\g$, as we will shortly see. 

\nl
We collect in the following two elementary lemmata what is needed in evaluating Lipschitz constants in Kolmogorov's scheme. 
Let $f$ be real analytic on $W_{\x,\e_0}$ depend also on $\o\in \O\subset \real^d_{\g,\x}$ and  assume it is  uniformly Lipschitz in $\o$, i.e.:
\beqno
\Lip_{\x,\e_0}(f):=\sup\frac{|f(y,x,\o)-f(y,x,\o')|}{|\o-\o'|}<\io\,,
\eeqno
where the supremum is taken over all $\o\neq\o'\in \O$ and over all $(y,x,\e)\in W_{\x,\e_0}$.

\lem{lip.1} Let $f$ as above, let  $\l= \Lip_{\x,\e_0}(f)$, and  let $0<\d<\x$. Then, the following holds.
\\
{\rm (i)} Let\footnote{$\natural_0=\{0,1,2,3,...\}$} $\a,\b\in \natural_0^d$ be multi--indices. Then, 
\beqno
 \Lip_{\x-\d,\e_0}(\partial_y^\a \partial_x^\b f)\le c\, \d^{-(|\a|+|\b|} \l\,,
 \eeqno
 for a suitable constant depending only on $d$ and $|\a|+|\b|$.
 
 \nl
{\rm (ii)} $\forall n\in\integer^d$, the  Fourier coefficients $f_n(y,\o)$ of $x\mapsto f(y,x,\o)$ satisfy\footnote{As usual, in Fourier analysis, the norm in the exponents are 1--norms.}
\beq{lem.2}
|f_n(y,\o)-f_n(y,\o')|\le \l e^{-|n|\x} \ |\o-\o'|\,,\quad \forall\ \o,\o'\in \O\,.
\eeq
{\rm (iii)}   
Assume $f_0(y,\o)=\langle f(y,\cdot,\o)\rangle=0$, for all $(y,x)\in D^d_\x\times\O$. 
Then, $F(y,x,\o):=D_\o^{-1}  f(y,x,\o)$ is Lipschitz in $\o$ and
\beqno
|F(y,x,\o)-F(y,x,\o')|\le \l' |\o-\o'|\,,\quad \forall\ y\in D^d_\x\,,\ x\in \torus^d_{\x-\d}\,,\ \o,\o'\in \O\,,
\eeqno
where, for   suitable constants\footnote{We take $k\ge k_1$ where $k_p$ is as in the `small--divisor' estimate Eq~(6) of \cite{C08}.} $c,k$  depending only on $d$ and $\t$, 
\beq{lem.4}
\l':= c\, \d^{-k}\g^{-2} (m+\l\g)\,,\quad m:=\sup_{W_{\x,\e_0}\times\O}|f|\,.
\eeq
\elem

\proof
(i) follows immediately by standard Cauchy estimates.

\nl
(ii) follows immediately by the standard ($n$--dependent) shift--of--contour argument based on Cauchy theorem of complex analysis, observing that
$$
f_n(y,\o)-f_n(y,\o')=\igl{\torus^d}{} \Big( f(y,x,\o)-f(y,x,\o')\Big)e^{-i n\cdot x} \frac{dx}{(2\p)^d}\,.
$$
and that $| f(y,x,\o)-f(y,x,\o')|\le \l |\o-\o'|$ on $W_{\x,\e_0}\times \O$.

\nl
(iii)  By \equ{Dom}, \equ{lem.2}, one has,  $\forall\ y\in D^d_\x$, $x\in \torus^d_{\x-\d}$, and $\o,\o'\in \O$:
\beqano
|F(y,x,\o)-F(y,x,\o')|&= &
\Big|
\sum_{n\neq 0}\Big( f_n(y,\o) \frac{(\o'-\o)\cdot n}{(\o\cdot n)(\o'\cdot n)}+ \frac{ f_n(y,\o)-f_n(y,\o')}{\o'\cdot n}
\Big) e^{i n\cdot x}\Big|\\
&\le&|\o-\o'|\sum_{n\neq 0}\Big( me^{-|n|\x}\frac{|n|^{2\t+1}}{\g^2}+\l\frac{|n|^\t}{\g}e^{-|n|\x}\Big) e^{|n|(\x-\d)}\,,
\eeqano 
which, since  $\o\in \real^d_{\g,\t}$, yields the claim by standard estimates\footnote{See, e.g., footnote 10 in \cite{C08}.}. \qed

\lem{lip.2} {\rm (i)} Let  $\O\subset \real^d$ and let $A=A(\o)$ be an invertible  matrix such that 
$\|A(\o)-A(\o')\|\le \l |\o-\o'|$ and  $\|A^{-1}(\o)\|\le m$, $\forall$ $\o,\o' \in \O$. Then, 
\beqno
\|A^{-1}(\o)-A^{-1}(\o')\|\le \l' |\o-\o'|\,,\quad \forall\ \o,\o'\in\O\,,
\eeqno
with $\l'=\l m^2$.

\nl
{\rm (ii)} For any $\o\in \O\subset \real^d$ and any $|\e|<\e_*$, let $ x\in\torus^d\mapsto\f(x)=\f(x,\o)=x+\e a(x,\o)\in\torus^d$ be  a near--to--identity $C^1$  diffeomorphism\footnote{The dependence upon $\e$ of the functions is not explicitly indicated.},  with inverse  given by $\psi(x')=\psi(x',\o) = x'+\e \a(x',\o)$, and satisfying  $\e_* \|a_x\|_{\!{}_\io}<1$. Assume that $|a(x,\o)-a(x,\o')|\le \l |\o-\o'|$ for any $x,\o,\o'$. Then, 
\beqno
|\a(x',\o)-\a(x',\o')|\le \l' |\o-\o'|\,,\qquad \forall \ \o,\o'\in \O\,,
\eeqno
with $\l'=\l (1-|\e| \|a_x\|_{\!{}_\io})$. Analogous statement holds in complex neighbourhoods of $\torus^d$.
\elem 

\proof (i) Let $v\neq 0$ and let $v'=A^{-1}(\o) v$. Then, for any $\o,\o'\in \O$, 
\beqano
|A^{-1}(\o) v- A^{-1}(\o') v| &=&
\big|A^{-1}(\o') \, \big( A(\o')  - A(\o)\big) v'\big|
\le m^2 \l |\o-\o'| \, |v|\,.
\eeqano
(ii) Let $x_1=\psi(x'_1,\o_1)$ and $x_2=\psi(x_1',\o_2)$. Then,
\beqano
|\a(x_1',\o_1)- \a(x_1',\o_2)|&=&|a(x_1,\o_1)- a(x_2,\o_2)|\\
&=& |a(x_1,\o_1)- a(x_1,\o_2) + a(x_1,\o_2)- a(x_2,\o_2)|\\
& \le& \l |\o_1-\o_2|+ \|a_x\|_{\!{}_\io} |x_1-x_2|\\
&=& \l |\o_1-\o_2|+ |\e| \|a_x\|_{\!{}_\io} |\a(x_1',\o_1)-\a(x_1',\o_2)|\,,
\eeqano
which implies the claim. The complex case is treated in the same way. \qed

\nl
To describe the iterative step needed to control Lipschitz constant in Kolmogorov's scheme we refer to \cite{C08} and, in particular, to 
Lemma 4 and\footnote{There is a small correction to be done in the statement of Lemma~4 in \cite{C08}, namely, the bound on $\|P'\|_{\bar\xi}$ in Eq.~(18) should be given {\sl after} hypothesis (19).}  its proof in \cite{C08}.

\pro{pro-lip} Let $E$, $Q$,  $T$ and $P$ be as in {\bf (2.2)} above, and assume that they  
depend in a Lipschitz way on $\o\in \O_\g$ with uniform (on their complex domain of definition) Lipschitz constant $\L$.
Let
\beqno
C:=\max\big\{|E|, |\o|,\|Q\|_{\x,\e_0}, \|T\|\,,\L\,,1\big\}\,,
\eeqno
assume that\footnote{This assumption, which is armless (since, eventually, $\g$ will be chosen small with $\e$), is made to simplify the estimate \equ{lem.4}.} $\g\le 1/2 \min\{1,\L\}$, and  let $0<\d<\x<1$. Finally, let $L$, 
$\phi_1={\rm id}+\e \tilde\phi$, $E_1=E+\e \tilde E$, $Q_1=Q+\e \tilde Q$, $T_1= T+ \e \tilde T$, and $P_1$ be as in step~(i) and Lemma~4 of {\rm \cite{C08}}. Then,  $\tilde E$, $\tilde Q$, $\tilde T$ and $\tilde \phi$ are  Lipschitz in $\o\in \O_\g$ uniformly on  $W_{\bar \x,\e_0}$, $\bar \x:=\x -\frac23 \d$, with Lipschitz constant given by   
\beqno
\tilde \L=c' \g^{-a'} C^{\m'} \d^{-\nu'} M\ge L\,,\qquad M:=\sup_{W_{\x,\e_0}\times \O_\g}\|P\|\,,
\eeqno
where 
$c', a', \m'\, \nu'$ are suitable positive constants  depending  on $\t$, $d$. Furthermore, if $\e_*\le \e_0$ is such that $\e_* \tilde \L\le \d/3$, then $P_1$ is 
Lipschitz in $\o\in \O_\g$ uniformly on $W_{\x',\e_*}$ with $\x':=\x-\d$ with Lipschitz constant $\tilde \L M$.
Finally, for any $|\e|< \e_*$,  $E_1$, $Q_1$ and $T_1$ are uniformly Lipschitz in $\o\in \O_\g$ with Lipschitz constant $\L_1:= \L+|\e| \tilde \L$.
\epro
\proof
Let us  give the details for the estimate on the Lipschitz constant of $\tilde E$.

\nl
$\tilde E$ is defined as\footnote{Compare step (i) at p. 135 of \cite{C08} and recall that $T(\o):=\langle Q_{yy}(0,\cdot;\o)\rangle^{-1}$.
Usually, we do not indicate the dependence upon $\e$ of the various functions involved.} $\o\cdot b+ P_0(0;\o)$
where 
$$
b=- T(\o) \big( \langle Q_{yy}(0,\cdot;\o) s_x \rangle +\langle P_y(0,\cdot;\o)\rangle\big)\,,\qquad
s(x;\o):= -D_\o^{-1} \big( P_y(0,\cdot;\o)-P_0(0;\o)\big) \,.
$$
Then, by Lemma~\ref{lip.1}--(iii) with $f=P_y(0,\cdot;\o)-P_0(0;\o)$, $\l=\L$, and using that $\L \g<1\le M$, we get\footnote{We denote possibly different constants depending on $d$ and $\t$ by $c$.} 
\beqno
\Lip_{\x-\frac\d3,\e_0}(s)\le c\, \d^{-k}\g^{-2} M\,,
\eeqno
and by  Lemma~\ref{lip.1}--(i), 
\beqno
\Lip_{\x-\frac{2\d}3,\e_0}(s_x)\le c\, \d^{-(k+1)}\g^{-2} M\,.
\eeqno
Now, by Lemma~\ref{lip.2}--(i) we get $\Lip(T)\le C^3$, and therefore\footnote{For products, use $\Lip(fg)\le \Lip(f) \sup|g|+ \Lip(g)\sup|f|$, and observe that $\|s_x\|_{\x-2/3\d}\le c\, \d^{-b}\g^{-1} M$; compare, e.g., Eq.~(6) in \cite{C08}.}
\beqno
\Lip_{\x-\frac{2\d}3,\e_0}(b)\le c\, C^4 \d^{k+1} \g^{-2} M\,, \qquad \Lip_{\x-\frac{2\d}3,\e_0}(E_1)\le  \L+ |\e| \cdot (c\, C^5 \d^{k+1} \g^{-2} M)\,.
\eeqno
It is not difficult to check that also the Lipschitz constants of\footnote{Notice that, since  $\tilde \L\ge L$, the hypotheses of Lemma~4  (compare Eq.~(19)) are met.}
 $\tilde Q$, $\b_0$, $\b$ (defined in Remark~2 (a) of \cite{C08}),  $\tilde T$ and $\tilde \phi$ satisfy similar estimates; also the estimate on 
$\Lip(P_1)$ is of the same type, but with an extra factor $M$, since in the definition of $P_1$ there appears  a term ($P^{(1)}$ in Eq.~(11) in \cite{C08}), which is quadratic in $\b$. \qed

\nl
Now, the inductive argument  follows  easily as in the proof of  Lemma~5 of \cite{C08}. We give a sketch of it.

\nl
Let, as in Lemma 5 of \cite{C08}, $\x_{j+1}=\x_j-\d_j$, $\d_j=\d_0/2^j$, $\d_0=(\x-\x_*)/2$, 
\beqno
\phi_j:W_{\x_j,\e_*}\times \O_\g\to D^d_{\x_{j-1}}\times \torus_{\x_{j-1}}^d\,,
\quad
\Phi_j=\Phi_{j-1}\circ \phi_j\,,\qquad (j\ge 1\,, \,\Phi_0= {\rm id})\,,
\eeqno 
so that $\phi_*$ in {\bf (2.2)} above is given by $\phi_*=\lim \Phi_j$. From the proof of Lemma~5 in \cite{C08} and from Cauchy estimate it follows that 
\beq{nab}
\sup_{B\times\torus\times \O_\g} \|\partial_z \Phi_j\|\le 2\,,\qquad z=(y,x)\,.
\eeq
Let 
\beq{liptil}
\tilde \L_i=c' \g^{-a'} C^{\m'} \d_i^{-\nu'} M_i\ge \Lip_{\x_i,\e_*}(\phi_i)\,,
\eeq
be as in the $i^{\rm th}$ iteration of\footnote{For $i\ge 0$,  
$E$, $\tilde E$, $Q$, $\tilde Q$, $\L$, $\tilde \L$,...,$\x$, $\d$, $\e$  correspond to 
$E_i$, $\tilde E_i$, $Q_i$, $\tilde Q_i$, $\L_i$, $\tilde \L_i$,...,$\x_i$, $\d_i$, $\e^{2^i}$, while 
$E_1$, $Q_1$, $\L_1$,... correspond to  $E_{i+1}$, $Q_{i+1}$, $\L_{i+1}$, etc.
} Proposition~\ref{pro-lip}, and let
let $\l_j:=\Lip(\Phi_j)$ be the Lipschitz constant of $\Phi_j$ over $B\times \O_\g$.
\\
Let $\o,\o'\in\O_\g$,  $z_j=\phi_j(y,\theta;\o)$ and $z'_j=\phi_j(y,\theta;\o')$, for $y\in B$ and $\theta\in\torus^d$. 
  Then, by \equ{nab} and \equ{liptil}
 \beqano
 |\Phi_j(y,\theta;\o)- \Phi_j(y,\theta;\o')|&=&|\Phi_{j-1}(z_j;\o)- \Phi_{j-1}(z'_j;\o')|\\
 &\le&|\Phi_{j-1}(z_j;\o)- \Phi_{j-1}(z'_j;\o)|+|\Phi_{j-1}(z_j';\o)- \Phi_{j-1}(z'_j;\o')|\\
 &\le& 2 |\phi_j(z,\theta;\o)- \phi_j(z,\theta;\o')|+ \l_{j-1} |\o-\o'|\\
 &\le& 2 |\e|^{2^j}\tilde \L_j |\o-\o'| +\l_{j-1}|\o-\o'|\,,
 \eeqano
which (dividing by $|\o-\o'|$ and taking the supremum over $y\in B$, $\theta\in\torus^d$ and $\o\neq \o'$), yields the relation
$$
\l_j\le\l_{j-1}+  2 \e^{2^j}\tilde \L_j\,,$$
which, iterated, implies\footnote{The super--exponential series is treated as in \cite{C08} p. 138.}, for $|\e|$ small enough,
$$
\l_j\le 1+ 2\sum_{i=0}^\io \e^{2^i}\tilde \L_i< 2  \,,\qquad \forall j\,.
$$
Taking the limit as $j\to\io$, we get  $\Lip(\phi_*)\le 2$, which, 
as discussed above, is all what is needed to conclude the proof of Theorem~2. \qed

\appendix

\section{An interview to Ya. Sinai}

One of the authors (I.F.), during her doctoral thesis \cite{f}, supervised by Luca Biasco and Ana 
Mill\'an Gasca, had the opportunity to interview Yakov Sinai on May 28, 2021, in his quality of student and witness to Kolmogorov's legacy. Here, we provide the transcription of this interview\footnote{\textbf{F}  = Isabella Fascitiello;   \textbf{S} = Yakov Sinai; \textbf{B} = Luca Biasco.}.
\\[1em]
\noindent
\textbf{F}:
The first question concerns Siegel's work on Diophantine estimates. These techniques are also used by Kolmogorov in his proof of the theorem in 1954, but he did not mention Siegel in the bibliography. Do you know if Kolmogorov was aware of Siegel's work on such matter?
\\[1em]
\noindent
\textbf{S}:
In my opinion, he didn't know Siegel's work. Siegel's work was discussed later in Arnol'd's seminar, and I assume that Arnol'd explained Siegel's work to Kolmogorov. As you know, they both used small denominators.
\\[1em]
\noindent
\textbf{F}:
Do you know what inspired Kolmogorov for Diophantine estimates?
\\[1em]
\noindent
\textbf{S}:
I'm not so sure about this.
\\[1em]
\noindent
\textbf{F}:
Okay. So, I move on to the next question: In the published text of the Amsterdam conference, Kolmogorov cited in bibliography `Mathematische Grundlagen der Quantenmechanik' (1932) by Von Neumann. Did Kolmogorov ever work on problems in quantum mechanics?
\\[1em]
\noindent
\textbf{S}:
Kolmogorov never worked on problems of quantum mechanics because he used to say that he didn't find interesting problems for himself in that field.
\\[1em]
\noindent
\textbf{F}:
Okay, but I have a puzzle to solve. I read a sentence written by Kolmogorov that I quote here: ``My papers on classical mechanics appeared under the influence of von Neumann's paper on the spectral theory of dynamical systems..."\footnote{\cite{selworks}, p 521}. And also in this sentence the reference is `Mathematical Foundations of Quantum Mechanics'.
\\[1em]
\noindent
\textbf{S}:
No, I remember it was another Von Neumann's paper; there was a paper written by Von Neumann about the ergodic theory.
\\[1em]
\noindent
\textbf{F}:
Okay. Actually, in another note, written by Shiryaev\footnote{In \cite{s}, p.53.} on Kolmogorov, the author wrote that there is another reference, that is `Operator Methods in Classical Mechanics'\footnote{\textit{Zur Operatorenmethode In Der Klassischen Mechanik}. Princeton, Annals of Mathematics, Second Series 33(3), (1932) pp. 587--642.}.
\\[1em]
\noindent
\textbf{S}:
It's possible. That was the main contribution in operator method.
\\[1em]
\noindent
\textbf{F}:
So it is correct the operator methods, and it is wrong the quantum mechanics?
\\[1em]
\noindent
\textbf{S}:
Yes, I think so.
\\[1em]
\noindent
\textbf{B}:
Professor Sinai, do you think that Kolmogorov, for his theorem on the persistence of invariant tori, was  also motivated by the foundations of statistical mechanics?
\\[1em]
\noindent
\textbf{S}:
He never mentioned this. He just mentioned the work of Chazy\footnote{Jean--Fran\c{c}ois Chazy (1882--1955). Two Chazy's papers, 1929 and 1932, are included in the references of \cite{K57}, both titled \textit{Sur l'allure finale du mouvement dans le probl\`eme des trois corps}.}. Chazy was a friend, a mathematician, or maybe a physicist; he was the first person who wrote about statistical and central limit theory and other papers on probability theory, which can be used in classical mechanics.
\\[1em]
\noindent
\textbf{F}:
Another question concerns your article  in the book \textit{Kolmogorov in Perspective}. You wrote that ``in the fall of 1957 Kolmogorov began a famous course of lectures on the theory of dynamical systems" and that, I quote, ``Kolmogorov presented the theorem\footnote{He refers to the theorem on the persistence of invariant tori for quasi--integrable Hamiltonian systems, Theorem~1 in \cite{K}.} that was to become the basis for the famous KAM theory, together with a complete proof".\footnote{See the complete excerpt cited in the first section 1 of this article, taken from \cite{s}.} What did you mean? The history of science says that the first complete proof is due to Arnol'd in 1963.
\\[1em]
\noindent
\textbf{S}:
There is a very good proof of Kolmogorov's theorem given by a student of Gallavotti\footnote{It refers to Giovanni Gallavotti (1942 -- ).}. I forgot his name.
\\[1em]
\noindent
\textbf{B}:
Maybe it's Luigi Chierchia.
\\[1em]
\noindent
\textbf{S}:
Maybe it was him. Yes.
\\[1em]
\noindent
\noindent
\textbf{F}:
But in 1957, Kolmogorov did present a complete proof in this seminar? Is this assertion true?
\\[1em]
\noindent
\textbf{S}:
You see, there is a controversy about this. For example, Arnol'd thought that Kolmogorov did not give a complete proof, that his proof had some gaps. And this was a reason why Arnol'd wrote his paper.
\\[1em]
\noindent
\textbf{F}:
Okay. What about you? Do you think Kolmogorov give a complete proof of his theorem?
\\[1em]
\noindent
\textbf{S}:
It is a controversial question. I believed that Kolmogorov gave a complete proof, but Arnol'd convinced me that Kolmogorov's proof was not complete.
\\[1em]
\noindent
\textbf{B}:
According to Arnol'd, was the proof incomplete because Kolmogorov omitted certain steps, or were there indeed certain gaps that Kolmogorov did not address?
\\[1em]
\noindent
\textbf{S}:
This is a complicated matter. There were some gaps in the estimates of the measure of invariant sets. That was the main point where Arnol'd complained about the proof by Kolmogorov. In Kolmogorov paper, complete estimates of such 
a 
measure were not given.
\\[1em]
\noindent
\textbf{B}:
 So, only regarding this specific point?
\\[1em]
\noindent
\textbf{S}:
Yes.
\\[1em]
\noindent
\textbf{F}:
Another question, maybe the last, concerning the connection between Kolmogorov and Arnol'd. Did Arnol'd ever make a comparison between his form of the  theorem on the persistence of invariant tori and Kolmogorov's original one? What were his motivations for giving a different proof of this theorem?
\\[1em]
\noindent
\textbf{S}:
Arnol'd wrote a complete proof of Kolmogorov's theorem which was published in a Russian journal, and it was exactly motivated by the fact that the proof in Kolmogorov's paper was not complete.
\\[1em]
\noindent
\textbf{F}:
Was Arnol'd thinking about celestial mechanics?
\\[1em]
\noindent
\textbf{S}:
Arnol'd continued thinking about celestial mechanics, but there was another student of Kolmogorov diligently working on the subject. I am specifically referring to Sitnikov\footnote{It refers to Kirill Aleksandrovich Sitnikov (1926 -- ?). Sitnikov is cited by Kolmogorov in the conference text \cite{K57}, with the article ``On Possible Capture in the Three--Body Problem", published in the Russian Journal \textit{Matematicheskii Sbornik} in 1953.}, who, in one of his articles, provides a comprehensive example of the solution to oscillations.
\\[1em]
\noindent
\textbf{F}:
So, Kolmogorov also was thinking about celestial mechanics in 1954?
\\[1em]
\noindent
\textbf{S}:
He was very much interested in problems in celestial mechanics. Alekseev's papers\footnote{
It refers to Vladimir Mikhailovich Alekseev (1932--1980), a student of Kolmogorov specializing in celestial mechanics To cite a few articles on the matter: ``Quasirandom vibrations and the problem of capture in the bounded three--body problem" (1967) in \textit{Doklady Akademii Nauk SSSR}, ``On the possibility of capture in the three--body problem with a negative value for the total energy constant" (1969) in \textit{Uspekhi Matematicheskikh Nauk}, and `Final motions in the three--body problem and symbolic dynamics' (1981) in \textit{Russian Mathematical Surveys}.} on celestial mechanics was certainly influenced by discussions with Kolmogorov.

\vglue1.truecm
\noindent
\footnotesize{(References are listed in chronological order)}

\end{document}